\documentclass[11pt]{article}

\usepackage[hmargin=3cm,vmargin=2cm,includefoot,includehead]{geometry}

\usepackage{mathtools,amssymb}

\title{On the independence of Robinson's set of axioms\\for propositional calculus}
\author{Beno\^it Jubin}
\date{29 September 2021}

\begin{document}

\maketitle

\begin{abstract}
We give a normal five-valued truth table proving independence of one of the axioms in Robinson's set of axioms for propositional calculus from 1968, answering a question raised in his article, where he uses a non-normal table.
We also give a normal four-valued table proving independence of one of the other axioms, where he uses a normal five-valued table.
\end{abstract}

In his article~\cite{robinson}, Thacher Robinson introduces a set of axioms for propositional calculus and proves their independence using many-valued truth tables.
In one of his independence proofs (for the axiom~(\texttt{S}) below), he uses a non-normal four-valued table, and he then asks in the last paragraph for a normal table showing independence.\footnote{All these concepts are explained in his article.
Normality is a property of a truth table for implication which is equivalent to the associated logic satisfying modus ponens.
We may use the singular ``table'' to denote a set of tables, one for each connective, and say that this set is normal when its table for implication is normal.}
He writes that Paul Bernays constructed for him such a normal six-valued table that was subsequently lost, and that apparently no normal five-valued tables exist.
In this note, I give a normal five-valued table proving the required independence.
The fact that it is normal implies the existence of a five-valued logic satisfying modus ponens and all of Robinson's axioms except that one.

The question is that of proving the independence of~(\texttt{S}) in Robinson's set of axioms:\footnote{Robinson gives two sets of axioms, one containing negation and the other containing falsum, but we can join these axioms in one set and consider both negation and falsum as primitive, instead of defining one of these connectives from the other and implication.}\footnote{Robinson uses a commuted form of \texttt{orelim} but the tables given here satisfy both forms.}

\begin{align}
(\texttt{K}) &&& p \to q \to p\\
(\texttt{S}) &&& (p \to q \to r) \to (p \to q) \to p \to r\\
\texttt{peirce} &&& ((p \to q) \to p) \to p\\
\texttt{andelimr} &&& p \wedge q \to p\\
\texttt{andeliml} &&& p \wedge q \to q\\
\texttt{andintro} &&& p \to q \to p \wedge q\\
\texttt{orintror} &&& p \to p \vee q\\
\texttt{orintrol} &&& p \to q \vee p\\
\texttt{orelim} &&& p \vee q \to (p \to r) \to (q \to r) \to r\\
\texttt{contrap} &&& (p \to \neg q) \to q \to \neg p\\
\texttt{notelim} &&& \neg p \to p \to q\\
\texttt{falseelim} &&& \bot \to p
\end{align}
where the sole inference rule is modus ponens
\begin{equation}
\texttt{mp} \quad\qquad\qquad\qquad p \quad \& \quad p \to q \quad \Longrightarrow \quad q.
\end{equation}

An answer is given by the following normal five-valued table, where the only designated value (the value considered true) is~0:
\begin{center}
\begin{tabular}{|c|ccccc|}
\hline
$\to$&0&1&2&3&4\\
\hline
0&0&1&1&1&1\\
1&0&0&0&0&0\\
2&0&0&0&0&0\\
3&0&0&4&0&4\\
4&0&0&3&3&0\\
\hline
\end{tabular}
\quad
\begin{tabular}{|c|ccccc|}
\hline
$\wedge$&0&1&2&3&4\\
\hline
0&0&1&1&1&1\\
1&1&1&1&1&1\\
2&1&1&1&1&1\\
3&1&1&1&1&1\\
4&1&1&1&1&1\\
\hline
\end{tabular}
\quad
\begin{tabular}{|c|ccccc|}
\hline
$\vee$&0&1&2&3&4\\
\hline
0&0&0&0&0&0\\
1&0&1&1&1&1\\
2&0&1&1&1&1\\
3&0&1&1&1&1\\
4&0&1&1&1&1\\
\hline
\end{tabular}
\quad
\begin{tabular}{|c|c|}
\hline
&$\neg$\\
\hline
0&2\\
1&0\\
2&0\\
3&1\\
4&1\\
\hline
\end{tabular}
\quad
\begin{tabular}{|c|}
\hline
$\bot$\\
\hline
1\\
\hline
\end{tabular}
\end{center}
which validates modus ponens and all axioms except~(\texttt{S}), which is for instance false for the assignment
\begin{equation}
[p \leftarrow 3, q \leftarrow 0, r \leftarrow 2].
\end{equation}

To prove independence of~(\texttt{K}), Robinson uses a normal five-valued table.
I found a normal four-valued table, where the designated values are 0, 1, 2:
\begin{center}
\begin{tabular}{|c|cccc|}
\hline
$\to$&0&1&2&3\\
\hline
0&0&0&2&3\\
1&0&0&3&3\\
2&0&0&0&3\\
3&0&0&0&0\\
\hline
\end{tabular}
\quad
\begin{tabular}{|c|cccc|}
\hline
$\wedge$&0&1&2&3\\
\hline
0&0&0&0&3\\
1&0&0&0&3\\
2&0&0&0&3\\
3&3&3&3&3\\
\hline
\end{tabular}
\quad
\begin{tabular}{|c|cccc|}
\hline
$\vee$&0&1&2&3\\
\hline
0&0&0&0&0\\
1&0&0&0&0\\
2&0&0&0&0\\
3&0&0&0&3\\
\hline
\end{tabular}
\quad
\begin{tabular}{|c|c|}
\hline
&$\neg$\\
\hline
0&3\\
1&3\\
2&3\\
3&0\\
\hline
\end{tabular}
\quad
\begin{tabular}{|c|}
\hline
$\bot$\\
\hline
3\\
\hline
\end{tabular}
\end{center}
which validates modus ponens and all axioms except~(\texttt{K}), which is for instance false for the assignment
\begin{equation}
[p \leftarrow 2, q \leftarrow 1].
\end{equation}

\section*{Addendum: another truth table}

Testing my program on other examples, I found a shorter normal truth-table for implication for the following independence problem.
In~\cite[\S3]{meyer}, Meyer and Parks use a normal four-valued table (with one designated value), found by Soboci\'nski, to prove independence of B' from \{W, pon, X\} where
\begin{align}
\text{B'} &&& (p \to q) \to (q \to r) \to p \to r,\\
\text{W} &&& (p \to p \to q) \to p \to q,\\
\text{pon} &&& p \to (p \to q) \to q,\\
\text{X} &&& ((((p \to q) \to q) \to p) \to r) \to ((((q \to p) \to p) \to q) \to r) \to r.
\end{align}
Here is a normal three-valued table where the designated values are 0, 1:
\begin{center}
\begin{tabular}{|c|ccc|}
\hline
$\to$&0&1&2\\
\hline
0&0&0&2\\
1&0&2&2\\
2&0&0&0\\
\hline
\end{tabular}
\end{center}
which validates modus ponens and \{W, pon, X\} but falsifies B' with the assignment
\begin{equation}
[p \leftarrow 1, q \leftarrow 0, r \leftarrow 1].
\end{equation}

\paragraph{Acknowledgments}
I would like to thank Norman Megill for bringing this problem to my attention, Mario Carneiro for independently confirming the truth tables, and Jean-Baptiste Bianquis for useful advice on OCaml matters.

\vspace*{5mm}
\noindent
\parbox[t]{21em}
{\scriptsize
Beno{\^i}t Jubin\\
Sorbonne Universit{\'e}\\
Institut de Math{\'e}matiques de Jussieu\\
4 place Jussieu, 75252 Paris Cedex 05, France
}

\end{document}